%Paper: 9205003
%From: Vladimir Pestov <vova@kauri.vuw.ac.nz>
%Date: Tue, 19 May 1992 04:26:24 GMT

%Nonstandard hulls of Banach-Lie groups and algebras. -
%Vladimir G. Pestov - Research Report RP-92-86, Department of
%Mathematics, Victoria University of Wellington,
%May 1992, 12 pp. - AmS TeX 2.1.

\documentstyle{amsppt}
\magnification=\magstep 1
\TagsOnRight
\NoBlackBoxes

% macros
\def\norm #1{{\left\Vert\,#1\,\right\Vert}}
\def\g {{\frak g}}
\def\R {^\star {\Bbb R}}
\def\m #1{{\mu_{\frak #1}(0)}}
\def\f #1{{fin~{\frak #1}}}
\def\h #1{\hat{\frak #1}}
\def\N{^\star{\Bbb N}}
\def\exp{{\widehat{exp}}}
\def\e{{\epsilon}}
% end of macros

\topmatter

\title
Nonstandard hulls of Banach-Lie groups and algebras
\endtitle

\author Vladimir G. Pestov
\endauthor

\affil
Department of Mathematics\\
Victoria University of Wellington \\
P.O. Box 600, Wellington, New Zealand
\endaffil

\abstract{We propose a new construction of
Banach-Lie groups and algebras via nonstandard analysis.
A major ``standard'' application is the Local Theorem which to
certain extent reduces the problem of associating a Lie group to a
given Banach-Lie algebra to a similar problem for finitely
generated Lie subalgebras.}
\endabstract
\endsubjclass

\endtopmatter
\document

\smallpagebreak
\heading  Introduction
\endheading
\smallpagebreak

The idea of forming the nonstandard hull of an internal normed space
(which is a cornerstone of our construction) belongs to nonstandard
Banach space theory. The notion, which appears nameless in the treatise of
 A. Robinson \cite{R1}
(Section 4.4), has received a comprehensive treatment
in the paper \cite{Lux1} by W.A.J. Luxemburg;
for a contemporary presentation of this area of functional analysis
and an extensive bibliography, the reader should consult \cite{HM}.

The procedure of forming the nonstandard hull may be also applied to
internal normed algebras. (And this has been done ---
in a different language, though --- in \cite{DCC}.)
In particular, the construction makes sense for an internal
normed Lie algebra. A less trivial observation is that a similar
construction is sensible relative to an internal Banach-Lie group, $G$,
an internal norm being fixed on its Lie algebra, $Lie~G$.
Under certain restricitons, the nonstandard hull, $\hat G$, of a
Banach-Lie group, $G$, becomes a Banach-Lie group in such a way that
its Lie algebra, $Lie~\hat G$, is isomorphic to the nonstandard
hull, $\widehat{Lie~G}$, of the internal Banach-Lie algebra
$Lie~G$.

One of the key differences between finite and infinite dimensional Lie
theories is an observation
that not every infinite dimensional Lie algebra is
isomorphic to the Lie algebra of a suitable Lie group.
Those Lie algebras isomorphic to Lie algebras of Lie groups
are called ``enlargable'' \cite{vEK, KYMO} or, less commonly, also
``globalizable'' \cite{Po} or
``integrable'' \cite{Le} (remark that the term
``integrable Lie algebra'' in \cite{Ka} has a somewhat
 different meaning).
Examples of non-enlargable Banach-Lie algebras may be found in
\cite{vEK, LaT, PrS}.

In this paper, the techniques of nonstandard hulls leads to
 new sufficient
conditions for the enlargability of a Banach-Lie algebra.
At the same time it turns out that the nonstandard hull of an
enlargable internal Banach-Lie algebra may be non-enlargable.

The paper is concluded with a discussion of possible applications and
further modifications of the construction.

I consider the present work as a modest contribution to
the investigations motivated by the Seventh Metamathematical Problem
of Abraham Rosinson \cite{R2}, and explain my viewpoint in more detail in
\cite{Pe3}.

The present investigation has been basically performed
as long ago as in
1987. In the same year a Russian version of this paper
 was submitted to and
by June 1988 accepted by the Siberian
Mathematical Journal; about the same time,
most of the main results have been announced in my brief note in French
\cite{Pe1}. For some reasons, later I have withdrawn the manuscript
 from the Sib. Math. J., and it is only
now that it appears in a published form.
 Some of the proofs in
my later papers \cite{Pe2, Pe7} were relying
on the results of the present work, and I regret  this
 delay in publication.

The results have been discussed at seminars at the Laboratory of
Functional Analysis (Institute of
Mathematics, Novosibirsk Science Center, Russian Federation,
then USSR; 1988) and
Laboratoire de Math\'ematiques (Universit\'e de Haute Alsace,
Mulhouse and Colmar, France; 1990),
and I am grateful to A.G. Kusraev and S.S. Kutateladze,
as well as
M. Goze, R. Lutz,  and A. Makhlouf, for their
kind hospitality at the corresponding
institutions. My thanks also go to I.V. L'vov and L.A. Bokut' who made a
 stimulating remark
after my talk at the 2nd Siberian School ``Algebra and Analysis''
(Tomsk; 1988).

In the paper I am working within the ``classical''  Robinson's
approach to nonstandard analysis \cite{R1, D, HL, LiR, LuG, SL}.
All algebras, groups, $\dots$ considered are assumed to belong
to a fixed set-theoretic structure $\frak M$. (I avoid using the term
``superstructure'' altogether, because at present the prefix ``super-''
definitely associates itself with a completely different sort of
mathematical objects.)
Unless otherwise stated, a nonstandard enlargement $^\star\frak M$
is supposed to be at least sequentially comprehensive
\cite{LiR} (for example, $\aleph_1$-saturated).

\smallpagebreak
\heading
1. Nonstandard hulls of normed Lie algebras
\endheading
\smallpagebreak

\subheading{(1.1)}
Let $\frak g$ be an internal normed Lie algebra.
Fix an internal norm, $\Vert\cdot\Vert$, on
$\g$  such that for a positive finite $M\in\R$,
the inequality $\norm{[x,y]}\leq M\cdot\norm x\cdot\norm y$
holds whenever $x,y\in\g$.
In particular, the submultiplicative (or ``admissible'') norms on Lie
algebras are those with $M=1$, while for the ``natural'' norms on
Lie algebras of operators on Banach spaces usually $M=2$
\cite{dlH}.

\subheading{(1.2)}
Denote
$$fin~\g =_{def} \{ x\in\g : \norm x ~{is ~a~finite~ element ~of}~\R\},$$
$$\mu_\g (0) =_{def} \{ x\in\g : \norm x \approx 0\}.$$
It is easy to verify that, by virtue of the condition (1.1),
 the {\it principal galaxy},
$fin~\g$, of a Lie algebra $\g$ is a (generally speaking, external)
Lie algebra, and the {\it monad} (or {\it halo}) of zero,
$\m g$, is a Lie ideal in $\f g$.

\subheading{(1.3)}
Denote by $\h g$ the quotient Lie algebra $\f g/\m g$
and call it the {\it nonstandard hull of the Lie algebra $\g$}
({\it with respect to the norm} $\norm\cdot$).
The corresponding quotient homomorphism will be denoted by
$\pi_\g$ or simply by $\pi : \f g \to \h g$.

The Lie algebra $\h g$ is endowed with a (standard $\Bbb R$-valued)
norm by letting $\norm y =_{def} st~\norm x$ for any element
$y\in\h g$ of the form $y=\pi x,~x\in \f g$. (The symbol $st$ stands
for the standard part map $fin~\R\to \Bbb R$.)

This definition is clearly consistent and it makes $\h g$ into a
standard normed Lie algebra with the following property:
for each $x,y\in\h g$, the inequality
$\norm{[x,y]}\leq st~M\cdot\norm x\cdot\norm y$ holds.

Since the underlying normed space of the Lie algebra $\h g$ is the
nonstandard hull of the  the underlying normed space of the Lie algebra
$\g$, then our assumption of the enlargement
$^\star\frak M$ being sequentially
comprehensive yields the completeness of $\h g$, that is,
$\h g$ becomes a Banach-Lie algebra. (See \cite{R1, HL, HM, Lux1, SL}.)

\subheading{(1.4)}
Denote by $U\equiv U_\g$ an open ball of the radius
$(1/3M)ln(3/2)$  in an internal complete normed
Lie algebra $\g$ centered at zero, where $M$ is as in (1.1).
Let $\hat U=_{def}U_{\h g}$ be a similar ball of the radius
$(1/3st~M)ln(3/2)$  in the hull Lie algebra
$\h g$.
Remark that $\pi_{\g} U = cl~\hat U$ and
$\pi_\g^{-1}\hat U$ is an external subset of the ball $U$ (with a
non-empty $S$-interior).

The Hausdorff series $h(x,y) = \sum_{r,s\in\N}\tilde H_{r,s}(x,y)$
converges within $U$, making $U$ into an internal Lie group germ
(our terminology and notation here follow \cite{B}, ch. II, \S 7).
Analogously, the Hausdorff series of the Lie algebra $\h g$, which we
denote by $h'(x,y) = \sum_{r,s\in\Bbb N}\tilde H'_{r,s}(x,y)$,
converges in $\hat U$; it makes the latter set into a (standard)
Lie group germ.

The following result is a key to the construction of the nonstandard
hull of a Banach-Lie group.

\proclaim
{(1.5). Proposition}
For each $x,y\in\pi^{-1}_\g \hat U$:

$$h'(\pi_\g x, \pi_\g y) = \pi_\g h(x,y).$$
\endproclaim
\demo{Proof}
First of all, remark that for all standard $r,s\in\Bbb N$:
$$\pi_\g \tilde H_{r,s}(x,y) = \tilde H'_{r,s}(\pi_\g x, \pi_g y)$$
(indeed, $\pi_\g$ is a Lie algebra homomorphism; now use a
manifest form of the Hausdorff series' components, as in
\cite{B}, ch. II, \S 6, $n^o$4, Th. 2).
Fix a $z\in U$ such that
$\pi_\g z = h'(\pi_g x, \pi_\g y)$.
For any $n\in\Bbb N$ standard, there exist numbers
$r(n),~s(n)\in\Bbb N$
such that for all standard $r_1>r(n)$ and $s_1>s(n)$:

$$\norm{\sum_{r=1}^{r_1}\sum_{s=1}^{s_1}\tilde H'_{r,s}(\pi_\g x,\pi_\g y) =
h'(\pi_\g x,\pi_\g y)}_{\h g}<{1\over n} \tag{1}$$

{}From (1) and the definition of the norm on $\h g$ it follows that for all
standard  $r_1>r(n)$ and $s_1>s(n)$:

$$\norm{\sum_{r=1}^{r_1}\sum_{s=1}^{s_1}\tilde H_{r,s}(x,y) - z}_{\h g}<
{1\over n} \tag{2}$$

Since the number sequence on the left hand side of (2) is internal,
then one may apply to it the First Principle of Permanence
(\cite{LiR}, ch. II, \S 7)
and conclude that there exist infinitely large positive integers
$\tilde r(n)$ and $\tilde s(n)$ such that (2) takes place whenever
$r(n)<r_1<\tilde r(n)$ and $s(n)<s_1<\tilde s(n)$.

The corollary 7.8 from the Principles of Permanence (\cite{LiR}, ch. II)
allows one to fix infinitely large positive integers $\tilde r$ and
$\tilde s$ with $\tilde r < \tilde r(n)$ and $\tilde s < \tilde s(n)$
for all $n\in\Bbb N$ standard.
For every $n\in\Bbb N$ standard the property (2) now holds if
$r_1=\tilde r$ and $s_1=\tilde s$. Therefore,

$$\norm{\sum_{r=1}^{r_1}\sum_{s=1}^{s_1}\tilde H_{r,s}(x,y) - z}_{\g}
 \approx 0, $$

\noindent that is,

$$\pi_\g (\sum_{r=1}^{r_1}\sum_{s=1}^{s_1}\tilde H_{r,s}(x,y)) =
h'(\pi_\g x, \pi_\g y) \tag{3}$$

We will now show that

$$\sum_{(r,s)\notin [0,\tilde r]\times [0,\tilde s]} \tilde H_{r,s}(x,y)
\approx 0 \tag{4} $$

Let us use the following majoration (\cite{B}, ch. II, \S 7, $n^o$2):

$$\norm{\tilde H_{r,s}(x,y)}\leq
\eta_{r,s}\cdot M^{r+s+1}\cdot\norm x\cdot\norm y ,\tag{5}$$

\noindent where $(\eta_{r,s})$ is a standard  convergent positive
number series, of which a particular form is of no interest to us.
By virtue of the sequence  $(\eta_{r,s})$ being standard, we have:

$$\sum_{(r,s)\notin [0,\tilde r]\times [0,\tilde s]} \eta_{r,s}
\approx 0 $$

\noindent (see, e.g., Th. 3.5.8 in \cite{D});
from the last property, the property (5), and the condition
$x,y\in U$ (which means that the quantity $\norm x\cdot\norm y$
on the right hand side of (5) is finite) the property (4) follows.

Finally, the desired equality is obtained by partitioning
the series $h(x,y)$ in the two in an obvious way and applying
to the two summands the properties (3) and (4).
\qed
\enddemo

\subheading{(1.6)}
The proposition (1.5) admits a clear heuristic reformulation of the
following kind: the nonstandard hull of a Lie group germ
associated to an internal Banach-Lie algebra, $\g$,
is naturally isomorphic to a Lie group germ associated
to the nonstandard hull of the internal Lie algebra $\g$
(``heuristic,'' because we do not define the notion of the
nonstandard hull of a Lie group germ).

\subheading{(1.7)}
As in the case of Banach spaces, each standard Banach-Lie algebra
$\g$ canonically embeds as a Banach-Lie subalgebra
into the nonstandard hull of its own, $\widehat{^\star\g}$.
This embedding is an isomorphism if and only if $\g$ is
finite-dimensional. In general, nonstandard hulls $\h g$ are
``very inseparable'' Banach-Lie algebras (cf. similar results
for Banach spaces \cite{HM}).

\smallpagebreak
\heading
2. Nonstandard hulls of Banach-Lie groups and enlargability
\endheading
\smallpagebreak

\subheading{(2.1)}
Let $G$ be an internal Banach-Lie group, let $\g$ be its Lie algebra,
and $exp : \g \to G$ be the coresponding exponential map. Fix on the
Lie algebra $\g$ an internal norm satisfying (1.1).

\subheading{(2.2)}
Let us call by the {\it principal galaxy} of the Lie group $G$
the set $fin~G = \cup \{ (exp~U_\g)^n : n\in\Bbb N\}$, and by the {\it monad}
 (or {\it halo}) {\it of the unit element} - the set
$\mu_G(e) =_{def} exp~\m g$.
(Recall that $U_\g$ is an open ball of the radius
$(1/3M)ln(3/2)$ centered at zero.)

\proclaim{(2.3). Proposition}
The principal galaxy of a Lie group $G$ is the intersection of all
(internal and external) subgroups of $G$ containing the image
of the principal galaxy of $\g$ under the exponential map.
\endproclaim
\demo{Proof}
Obviously, $fin~G$ is a semigroup which is contained in any subgroup
of $G$ with the above property.
It remains only to notice that $exp~U$ is a symmetric set
(because for $x\in U,~(exp~x)^{-1} = exp~(-x)$);
since the semigroup $fin~G$ shares this property then it is a group.
\enddemo

\proclaim{(2.4). Lemma}
Let $V$ be a ball about zero of a finite non infinitesimal radius  in $\g$.
Then the principal galaxy of the Lie group $G$ may be represented as
follows:

$$fin~G = \cup\{(exp~V)^n : n\in\Bbb N\}.$$
\endproclaim
\demo{Proof}
According to the condition (1.1), there exists an $n\in\Bbb N$ such that
$n^{-1}U\subset V \subset nU$.
Let $x\in fin~G$.
For some $x_1, x_2, \dots , x_k \in exp~U,~ x = x_1x_2\dots x_k$.
Pick $y_i\in U$ with $exp~y_i = x_i$ and set
$z_i = n^{-1}y_i$.
Clearly, $z_i\in V$ and $x = \prod^k_{i=1}(exp~z_i)^n$.
Thus, the inclusion $\subset$ holds. The inverse inclusion,
$\supset$, is proved similarly.
\qed
\enddemo

\subheading{(2.5)}
The binary operation in both Lie group germs $U$ and $\hat U$
will be denoted from now on by a dot $(.)$ as in \cite{B}.
Now the proposition (1.5) may be rewritten as follows:
for every $x,y\in\pi^{-1}_\g \hat U$,
$$ \pi_\g (x.y) =
(\pi_\g (x)).(\pi_\g (y))$$

\proclaim{(2.6). Lemma}
The monad of zero, $\m g$, of a normed Lie algebra $\g$ is a
subgroup of the Lie group germ $U$.
In addition, for any $x\in\pi^{-1}_\g (\hat U)$,
one has $x.\m g .(-x) \subset \m g$.
\endproclaim
\demo{Proof}
Let $x,y\in\m g$.
We use the following majoration (\cite{B}, ch. II, \S 7, $n^0$2):

$$\norm{x.y}\leq - {1 \over M} ln(2-exp~M(\norm x+\norm y)).$$

It follows that $\norm{x.y}\approx 0$, that is,
$x.y\in\m g$.
Now let $x\in\pi^{-1}_\g\hat U$ and let $y\in\m g$.
By virtue of (2.5), one has

$$\pi_\g [x.y.(-x)] = \pi_\g (x).\pi_\g (y).(-\pi_\g (x)) =
\pi_\g (x).0.(-\pi_\g (x)) = 0.$$

Therefore, $x.y.(-x)\in\pi^{-1}_\g (0) = \m g$.
\qed
\enddemo

\proclaim{(2.7). Proposition}
The monad of unit, $\mu_G(e)$, of a Banach-Lie group $G$ is a
normal subgroup of the principal galaxy $fin~G$.
\endproclaim
\demo{Proof}
It follows from (2.6) and the main properties of the
exponential map that $\mu_G(e)$ is a subgroup of $fin~G$ invariant
under all inner automorphisms made by elements of
$exp~(\pi^{-1}_\g(\hat U))$.
Since, for example, ${1\over 2}U\subset\pi^{-1}_\g(\hat U)$,
then the lemma (2.4) implies the normality of
$\mu_G(e)$ in $fin~G$.
\qed
\enddemo
\subheading{(2.8)}
We denote by $\hat G$ the quotient group $fin~G/\mu_G(e)$
and by $\pi_G$ the corresponding quotient homomorphism
$fin~G\to\hat G$.

\subheading{(2.9)}
Let us define a map $\exp~ : \h g\to \hat G$ by letting for each
$x\in \f g$:

$$\exp~\pi_\g x = \pi_Gexp~x.$$

\demo{Correctness of the definition}
If $\pi_g x =\pi_g y$ then for some $n\in\Bbb N,~ x_1 =_{def} n^{-1}x\in U$
and $y_1 =_{def} n^{-1}y\in U$.
Now one has $x_1.(-y_1)\in\m g$ (2.5) and further:
$$(\pi_Gexp~x)(\pi_Gexp~y)^{-1} =
(\pi_Gexp~x_1)^n(\pi_Gexp~y_1)^{-n} = $$
$$(\pi_Gexp~x_1)^{n-1}\pi_Gexp~(x_1.(-y_1))(\pi_Gexp~y_1)^{-n+1}=$$
$$(\pi_Gexp~x_1)^{n-1}(\pi_Gexp~y_1)^{-n+1}=\dots = e_{\hat G}.$$
\qed
\enddemo
\proclaim{(2.10). Lemma}
Let  $x\in\hat U,~y\in\hat U$,
and $x.(-y)\in\hat U$.
Then $\exp~ [x.(-y)] = \exp~ x$.
\endproclaim
\demo{Proof}
For arbitrary $x', y'\in U_\g$ such that $\pi_g x' = x$ and
$\pi_g y' = y$, the following takes place:

$$exp~[x'.(-y')] = exp~x'(exp~y')^{-1}.$$

It remains to apply the homomorphism $\pi_G$ on both sides of this
equality and then use (2.5) and (2.9).
\qed
\enddemo

\subheading{(2.11)}
Let us impose an additional condition upon the norm on a Lie algebra
$\g$ by requiring the existence of a standard $\e>0$
such that the restriction $exp_G\vert U_\e : U_\e\to G$ is one-to-one
(or, without a loss of generality, a local diffeomorphism).

\proclaim{(2.12). Lemma}
Under assumptions 2.1 and 2.11, the restriction of the mapping
$\exp$ to a neighbourhood of zero in $\h g$ is one-to-one.
\endproclaim
\demo{Proof}
Pick a neighbourhood $W$ of the origin in $\h g$ meeting the condition
$W.(-W)\subset \hat U\cap U_{{1\over 2}st~\e}$.
Suppose that for some $x_1, x_2\in W$ such that $x_1\neq x_2$,
one has $\exp ~x_1 = \exp ~x_2$.
Fix $y_1,y_2\in\pi_\g^{-1}W$ with
$\pi_\g y_i=x_i,~i=1,2$.
By virtue of (2.9) one concludes:

$$exp~[y_1.(-y_2)] = exp~y_1(exp~y_2)^{-1}\in\pi_G^{-1}(e)=\mu_G(e),$$

\noindent that is, for some $z\in\m g$, one has
$exp~[y_1.(-y_2)]=exp~z$.
At the same time, $y_1.(-y_2)\in U_\e$ and
$y_1.(-y_2)\notin \m g$ because $x_1\neq x_2$. This contradicts the
one-to-one property of $exp_G\vert U_\e$.
\qed
\enddemo

\proclaim{(2.13). Theorem}
Under the assumptions (2.1) and (2.11), there exists a unique structure
of a Banach-Lie group on $\hat G$ such that the mapping $\exp$ is
a local diffeomorphism. In this case the Lie algebra $Lie~\hat G$
is canonically isomoprhic to the Lie algebra $\h g$,
and $\exp$ coincides with the exponential map
$exp_{\hat G}$.
\endproclaim
\demo{Proof}
It suffices to apply to (2.10) and (2.12) the theorem on extension of
analytic structure \cite{\`S}.
The latter statement of the theorem follows from
(\cite{B}, ch. III, Th. 4, $n^o$(v)) and a definition of the exponential
map in the \S 4 {\it ibidem}.
\qed
\enddemo
\subheading{(2.14)}
We will call the Banach-Lie group $\hat G$ the {\it nonstandard hull
of an internal Banach-Lie group} $G$ formed w.r.t.
a fixed internal norm on the Lie algebra $\g \simeq Lie~G$.

\subheading{(2.15)}
How much restrictive is the condition (2.11), imposed on a norm on
a Lie algebra $\g$?

First of all, there is a classical result stating that if the connected
component of unit element, ${\Cal Z}_0(G)$, of the center, $\Cal Z(G)$,
of a Banach-Lie group $G$ is simply connected, that is,
$\pi_1\Cal Z_0(G)=(0)$, then for any norm on the Lie algebra
$\g\simeq Lie~G$, satisfying the condition (2.1), the condition (2.11)
follows automatically \cite{LaT}.

\subheading{(2.16)}
Otherwise, one should make sure that (2.11) be satisfied. A simplest example
is provided by the Lie group $^\star U(1)$.
Having realized the exponential map
$exp : ^\star\frak u(1) \simeq \R\to^\star U(1)$ as
$t\mapsto exp\{2\pi it\}$ and defining a norm on $^\star\frak u(1) \simeq \R$
by the rule $\norm t' =_{def}\alpha \vert t\vert$,
where $\alpha >0$ and
$\alpha\approx 0$, it is easy to show that

\item{a)} $\widehat{\R} \simeq \Bbb R$;
\item{b)} $fin~^\star U(1) = \mu_{^\star U(1)}(e) = ^\star U(1)$;

\noindent
therefore, the nonstandard hull $\widehat{^\star U(1)}$ formed w.r.t.
the above special norm $\norm\cdot '$ on $\R$
is trivial (coincides with $(e)$)
and the relation $Lie~\hat G \simeq \h g$ in this case cannot be achieved.

\subheading{(2.17)}
Now we aim at constructing an example (2.21) of an enlargable internal
Banach-Lie algebra $\g$ with a non-enlargable nonstandard hull, $\h g$.

First we recall a construction of a non-enlargable Banach-Lie algebra
due to Lazard and Tits \cite{LaT}.
Let $G$ be a connected simply connected Banach-Lie group such that
$\pi_1\Cal Z_0(G)\neq (0)$. (For example, the unitary group
$U(\Cal H)$ of an infinite-dimensional Hilbert space, $\Cal H$, with
the uniform operator topology has this property, see\cite{Ku}.)
There exists a one-dimensional toroidal subgroup $T<\Cal Z(G)$;
let $\frak t$ denote the corresponding Lie subalgebra of the center,
$\frak z$, of the Lie algebra $\g\simeq Lie~G$.
Denote by $\frak z_a$ a central ideal of the Lie algebra
$\g\oplus\g$ of the form $\frak z_a =_{def}\{(x, ax) : x\in \frak t\}$,
where $a\in\Bbb R$.
The heart of the Lazard---Tits construction is the fact that the quotient
Banach-Lie algebra $\g_a =_{def} (\g\oplus\g)/\frak z_a$ is
enlargable if and only if $a\in\Bbb R$ is a rational number.

\subheading{(2.18)}
The Lazard---Tits construction results in an observation that
if $\g$ is an enlargable Banach-Lie algebra such that every quotient Lie
algebra of $\g\oplus\g$ by a one-dimensional central subalgebra is
enlargable, then $\pi_1\Cal Z_0(G)= (0)$,
where $G$ stands for a connected
 simply connected Banach-Lie group associated
to $\g$.

\subheading{(2.19)}
In particular, it is well known that $\pi_1\Cal Z_0(G)= (0)$ for any
simply connected finite-dimensional Lie group $G$
(see \cite{LaT} or \cite{vEK}, property (0)).

\subheading{(2.20)}
If a Banach-Lie algebra $\g$ is solvable then it is enlargable according
to a theorem of \`Swierczkowski \cite{\`S}; by virtue of (2.18)
 this theorem
 implies that $\pi_1\Cal Z_0(G)= (0)$
for a connected simply connected Banach-Lie group $G$ attached to
a solvable Banach-Lie algebra $\g$.

\subheading{(2.21). Example}
Let $a$ be a finite number from $^\star\Bbb Q$ such that
$st~a\in\Bbb R$ is irrational.
Fix a norm meeting the condition (1.1)
on a standard Banach-Lie algebra $\g$ with the properties
 listed in (2.17).
(For example, for $\g = \frak u(\Cal H)$ it may be the operator norm
inherited from $\Cal{L(H)}$ with $M=2$.)
Extend this norm to $\g\oplus\g$ (say, as an $l_1$-type sum) and
endow the quotient Lie algebra $^\star\g_a=_{def}(^\star\g\oplus^\star\g)/
\frak z_a$ with the quotient norm.

What follows, are all easily verifiable statements.

\item{a)} The Lie algebra  $^\star\g_a$ is enlargable.
\item{b)} There exists a natural isometrical embedding
$j :\g^2\hookrightarrow\widehat{^\star\g^2}$
 (cf. (3.1) --- (3.2) in our paper below).
\item{c)} $j(\frak z_{st~a})=\widehat{\frak z_a}$.
\item{d)} For any normed algebra $\frak h$ and a Lie  ideal $\frak o$,
a canonical isometrical Lie algebra isomorphism takes place:
$\h h/\h o \simeq \widehat{\frak h/\frak o}$.
\item{e)} Therefore, the Lie algebra $\g_{st~a}$ isometrically embeds
into the nonstandard hull $\widehat{^\star\g_a}$.
\item{f)} To conclude that  $\widehat{^\star\g_a}$ is a
non-enlargable  Lie algebra, it suffices to remember that a closed Lie
subalgebra of an enlargable Banach-Lie algebra is enlargable
(\cite{vEK}, p. 22, item ($^{\star\star\star}$)).

\subheading{(2.22)}
At the same time, we will present two exclusions from a general rule
exemplified in (2.21).

Let us call a Banach-Lie algebra $\g$
{\it hyperfinite dimensional} if it is isomorphic to the nonstandard
hull of a $\star$finite dimensional Lie algebra. (This name agrees
with the terminology of nonstandard Banach space theory, cf. \cite{HM}.)

\proclaim{(2.13). Theorem}
Any hyperfinite dimensional Banach-Lie algebra $\g$ is enlargable.
\endproclaim
\demo{Proof}
It is sufficient to use the classical Lie-Cartan theorem together with
(2.19) and (2.13), applied to a $\star$finite dimensional Lie algebra
$\frak h$ such that $\h h\simeq \g$.
\qed
\enddemo

\subheading{(2.24)}
Let us call a Banach-Lie algebra $\g$ {\it hypersolvable}
if it is isomorphic to the nonstandard hull of an internal solvable
(= $\star$solvable) Banach-Lie algebra. In a manner similar to (2.23),
the following result is deduced from (2.20).

\proclaim{(2.25). Theorem}
Any hypersolvable Banach-Lie algebra is enlargable.
\endproclaim

\smallpagebreak
\heading
3. Applications to the enlargability of standard Banach-Lie algebras:
a local theorem
\endheading
\smallpagebreak

\subheading{(3.1)}
Let $\phi$ be an internal bounded homomoprhism from an internal normed
Lie algebra $\g$ to a similar Lie algebra $\frak h$; suppose that
there are norms fixed
on both Lie algebras $\g$ and $\frak h$ satisfying (1.1), and in addition
the operator norm $\norm\phi$ of the homomoprhism $\phi$ is a
finite element of $\R$.

Under these conditions, the homomoprhism $\phi$ determines a
homomoprhism between the nonstandard hulls,
$\hat\phi : \h g\to\h h$, by means of the rule
$\hat\phi (\pi_g x) = \pi_\frak h (\phi x)$ whenever
$x\in fin~\g$.
It is easy to verify that the definition is consistent and
that $\norm{\hat\phi}\leq st~\norm\phi$.

We will call $\hat\phi$ the {\it nonstandard hull of} $\phi$.

\subheading{(3.2)}
In particular, if $\phi : \g\to\frak h$
is an isometric embedding of normed Lie algebras then the nonstandard
hull $\hat\phi : \h g\to \h h$ is so.

\proclaim{(3.3). Lemma}
Let $\g$ be a standard normed Lie algebra and let $\frak h$ be an internal
normed Lie subalgebra of the Lie algebra $^\star\g$.
If for each $x\in\g$ the intersection $\mu_\g(x)\cap\frak h$ is non-empty
then the Lie algebra $\g$ embeds isometrically into the nonstandard hull
$\h h$.
\endproclaim
\demo{Proof}
Let $i : \g\to\widehat{^\star \g}$
be the canonical embedding (1.7). Denote by $j$ the embedding
$\frak h\hookrightarrow ^\star\g$;
then $\hat j : \h h\to \widehat{^\star \g}$ is an isometric embedding
(3.2). Now it follows from the conditions of our Lemma that for each
$x\in\g$ there exists an $y\in\frak h$ with
$i(x)=\pi_\frak h(y)$.
Therefore, $i(\g)\subset \hat j (\h h)$, that is,
$\g$ embeds isometrically into $\h h$ by means of the map
$j^{-1}i$.
\qed
\enddemo

Here is the main result of our paper.

\proclaim{(3.4). Local Theorem on Enlargability of Banach-Lie Algebras}
Let $\g$ be a Banach-Lie algebra.
Suppose that there exist a family $\frak H$ of closed Lie subalgebras
and a neighbourhood $V$ of zero such that:

\item{1)} For each $\frak h_1, \frak h_2\in\frak H$ there is
an $\frak h_3\in\frak H$ such that $\frak h_1\cup\frak h_2\subset
\frak h_3$;

\item{2)} $\cup\frak H$ is dense in $\g$;

\item{3)} every Lie algebra $\frak h\in\frak H$ is enlargable, and if
H is a corresponding connected simply connected Lie group then the
restriction $exp_H\vert V\cap\frak h$ is one-to-one.

Then the Lie algebra $\g$ is enlargable.
\endproclaim
\demo{Proof}
Let the nonstandard enlargement $^\star\frak M$ of an appropriate
set-theoretic structure $\frak M$ be $(Card~\g)^+$-saturated.

Pick a norm on the Lie algebra $\g$ satisfying (1.1); let
$\e\in\Bbb R_+$ be such that the $\e$-ball $U_\e(0)$ is contained in $V$.

The density of the set $X=_{def}\cup\frak H$ in $\g$ means the
non-emptiness of each intersection $^\star X\cap\mu_\g(x)$.
The theorem 7.2.6 in \cite{SL}
enables one to pick a $\star$finite subset $A\subset^\star X$
such that for each $x\in\g$, the intersection
$A\cap\mu_\g(x)\neq\emptyset$.

The condition 1) guarantees the existence of an internal Lie algebra
$\frak h\in^\star\frak H$
which contains $A$.
The $\star$Transfer and the condition 3) imply the enlargability of the
normed Lie algebra $\frak h$; moreover, the corresponding
internal Lie group, $H$,
may be chosen so as to make  $exp_H\vert ^\star V\cap\frak h$
one-to-one.

Clearly, the restriction of $exp_H$ to the $(\e/2)$-ball centered at
zero is one-to-one as well.
Therefore, $\frak h$ meets the condition (2.11), and  according to (2.13),
the nonstandard hull $\h h$ is an enlargable standard Banach-Lie algebra.

By virtue of (3.3), the Lie algebra $\g$ embeds in $\h h$ as a closed Lie
subalgebra; now it remains to use the hereditary property of
enlargability of Banach-Lie algebras w.r.t. closed Lie subalgebras
(\cite{vEK}, p. 22, item ($^{\star\star\star}$)) to state the
enlargability of $\g$ itself.
\qed
\enddemo

Taking into account (2.19) and (2.20), one ontains the following
corollaries.

\proclaim{(3.5). Corollary}
Suppose that a Banach-Lie algebra $\g$ contains a dense Lie subalgebra
$\frak h$ such that every finitely generated Lie subalgebra
of $\frak h$ is finite-dimensional.
Then $\g$ is enlargable.
\endproclaim
\qed

\proclaim{(3.6). Corollary}
Suppose that a Banach-Lie algebra $\g$ contains a dense Lie subalgebra
$\frak h$ such that every finitely generated Lie subalgebra
of $\frak h$ is solvable. (In other terms, $\g$ contains a dense
locally solvable subalgebra.)
Then $\g$ is enlargable.
\endproclaim
\qed

\smallpagebreak
\heading
4. Concluding remarks
\endheading
\smallpagebreak

\subheading{(4.1)}
At present, Banach-Lie groups and algebras do not play any
 noticeable role
within applied infinite-dimensional Lie theory.
(On an appropriate occasion, V.I. Arnol'd
\cite {A} has compared the Banach manifold
theory with a ``Procrustean bed.'')
Those infinite-dimensional Lie groups used in mathematical physics
either are not ``structurized'' at all or become Fr\'echet-Lie groups
\cite{KYMO} or even more general Lie groups modeled over arbitrary
locally convex spaces \cite{M}.

On the other hand, the theory of nonstandard hulls of Banach spaces
(or, just the same, the theory of ultraproducts of Banach spaces)
has been transferred to the case of general locally convex spaces
in recent years \cite{H}. In this connection, it should be of a certain
interest, to try to generalize the concept of the nonstandard hull to wider
classes of Lie groups.
Although a task constitutes no hardships for Lie algebras
by an analogy with \cite{H}, the case is different with Lie groups
because, for instance, the exponential map $exp_G$ may turn out
not to be a local diffeomorphism for a Fr\'echet-Lie group $G$
\cite{KYMO, M}.

\subheading{(4.2)}
I expect that the construction of nonstandard hull
 may find an application in the
quantum chromodynamics (QCD) while investigating the so-called
QCD-limit $N_c\to\infty$, where $N_c$ stands for the number of colors
of quarks. In the realistic QCD the number $N_c$ equals 3, and $SU(3)$ is
the color symmetry group of the charge space;
in the most general case, however,
the ``rigid'' gauge group is $SU(N_c)$.
Investigation of the QCD-limit is an important and still open
problem of the gauge theory (cf. the lecture 8 in \cite{N} and
remarks on this topic in \cite{Sch}).

The limit ``rigid'' gauge group $SU(\infty)$ is either treated
heuristically or assumed to coincide with the strict direct limit
$lim_{\to}SU(n)$, which is a Lie group modeled over an {\it (LB)}-space.
However, this Lie group, with its relatively poor structure,  can hardly
reflect all the complexity of the limit transition.

I propose to consider in this role the group $\widehat{SU(\nu)}$,
where $\nu$ is an infinitely large positive integer.
It is easily verifiable that $\widehat{SU(\nu)}$ acts on the
``charge space'' $\widehat{{\Bbb C}^\nu}$ which is an infinite-dimensional
Hilbert space.

The geometric formalism underlying the gauge theories (fiber bundles,
connexions etc.) seems to be readily amenable to the procedure of
forming nonstandard hulls.
Of course, one would arrive at ``global,'' infinite-dimensional
standard objects modeled on Banach spaces, as in the approach of
\cite{BdK}. Surely, the Loeb measure will play a central role
in the construction of the action functional.

Such an interpretation of the QCD-limit would be in a full accordance
with the principles of hyper quantum mechanics \cite{Ke}.

\subheading{(4.3)}
A nonarchimedean version of the concept of nonstandard hull, that of a
$\rho$-nonarchimedean hull, due to W.A.J. Luxemburg \cite{Lux2},
has also been extended to Banach-Lie algebras and groups \cite{Pe4, Pe6}.
Hopefully, this notion has something to do with deformation and
quantization of algebras and groups.

\subheading{(4.4)}
Finally, the techniques of nonstandard hulls
of Lie groups,
Lie superalgebras, and Grassmann algebras seems to be applicable
to certain problems
of analysis and geometry with anticommuting
variables (known as
superanalysis and supergeometry), see \cite{Pe2, Pe5, Pe7}.

\vfill
\eject

\enddocument
\bye